\documentclass[11pt,twoside]{article}
\usepackage{graphicx}
\usepackage{amssymb}
\usepackage{amsmath}
\usepackage{graphics}

\title{\bf {Application of optimal homotopy asymptotic method to nonlinear
Bingham fluid dampers}}

\author{Vasile Marinca$^1$, Remus-Daniel Ene$^2$, Liviu Bereteu$^3$}

\date{}

\begin{document}

\thispagestyle{plain}

\maketitle
\begin{center}
{$^1$ University Politehnica Timi\c soara, Department of Mechanics and Vibration, Timi\c soara, 300222, Romania\\
Department of Electromechanics and Vibration, Center for Advanced
and Fundamental Technical Research, Romania Academy, Timi\c soara,
300223, Romania, e-mail: vmarinca@mec.upt.ro} 
\end{center}

\begin{center}
{$^{2}$ University Politehnica Timi\c soara, Department of
Mathematics, Timi\c soara, 300006, Romania, e-mail:
remus.ene@upt.ro}
\end{center}

\begin{center}{$^3$ University Politehnica Timi\c soara,
Department of Mechanics and Strength of Materials, 300222 Timi\c
soara, Romania, e-mail: liviu.bereteu@upt.ro}
\end{center}

\begin{abstract}
Magnetorheological fluids (MR) are stable suspensions of
magnetizable microparticles, characterized by the property to
change the rheological characteristics when subjected to the
action of magnetic field. Together with another class of materials
that change their rheological characteristics in the presence of
an electric field, called electrorheological materials are known
in the literature as the smart materials or controlled materials.
In the absence of a magnetic field the particles in MR fluid are
dispersed in the base fluid and its flow through the apertures is
behaves as a Newtonian fluid having a constant shear stress. When
the magnetic field is applying a MR fluid behavior change, and
behaves like a Bingham fluid with a variable shear stress. Dynamic
response time is an important characteristic for determining the
performance of MR dampers in practical civil engineering
applications. The purpose of this paper is to show how to use the
Optimal Homotopy Asymptotic Method (OHAM) to solve the nonlinear
differential equation of a modified Bingham model with non-viscous
exponential damping. Our procedure does not depend upon small
parameters and provides us with a convenient way to optimally
control the convergence of the approximate solutions. OHAM is very
efficient in practice ensuring a very rapid convergence of the
solution after only one iteration and with a small number of step.

\end{abstract}
\bigskip

\section{Introduction}
\label{sec:1}

\noindent \par Over the past decade, much attention has been given
to the MR dampers for its attractive characteristics in
applications of civil engineering structures including earth-quake
hazard mitigation, or high strength, insensitivity to
contamination, and small power requirement. Also MR fluids have
attracted considerable interest due to their wide range of use in
vibration dampers for vehicle suspension systems or machinery
mounts, their stiffness and damping, characteristics can be
adjusted very quickly by applying a suitable electric or magnetic
field \cite{1}-\cite{4}. The magnetorheological response of MR
fluids results from the polarization induced in the suspended
particles by application of an external field. The interaction
between the resulting induced dipoles causes the particles to form
columnar structures, parallel to the applied field. These
chain-like structures restrict the motion of the fluid, thereby
increasing the viscous characteristics of the suspension. The
mechanical energy needed to yield these chain-like structures
increases as the applied field increases resulting in a field
dependent yield stress \cite{3}. Dynamic constitutive relation of
MR fluids is very complicated and provided damping force is
intrinsically nonlinear, so there is not a consistent recognized
mechanical model for MR dampers. The mechanical model for an MR
damper is often established through optimization method according
to experimental data \cite{5}-\cite{8}. At present, there is a
variety of dynamic models for MR dampers. Some models which are
simple cannot effectively simulate nonlinear dynamic
characteristics of MR dampers. Although there are some models
which can simulate nonlinear dynamic characteristics, they are
established by strong nonlinear equations having a lot of
parameters which result in complicated numerical calculation. In
civil engineering Bingham model are often used for emulating the
dynamic behavior of MR dampers. This model is one of the most
popular models have been widely used with reasonable accuracy and
computational cost. The so-called Bingham model includes a
variable rigid perfectly plastic element connected in parallel to
a Newtonian viscosity element. This model assumes that the fluid
exhibits shear stress proportional to shear rate in the post-yield
region and can be expressed as \cite{4}, \cite{6}, \cite{7},
\cite{9}:
\begin{eqnarray}\label{Rheo1}
\displaystyle   \tau = \tau_y(H) \cdot sgn(\dot{\gamma}) + \eta
\dot{\gamma}
\end{eqnarray}
where $\tau$ is the shear stress in the fluid, $\tau_y$ is the
yielding shear stress controlled by the applied field $H$, $\eta$
is the post-yield viscosity independent of the applied magnetic
field, $\gamma$ is the shear strain rate and $sgn( \cdot )$ is the
signum function. That is, the fluid is in a state of rest and
behaves viscoelastically until the shear stress is greater than
the critical value $\tau_y$, whereas it moves like a Newtonian
fluid when such a critical value is exceed.

Based on this model of the rheological behavior of
electrorheological fluids, Stanway et al. \cite{5} proposed an
idealized mechanical system which consists of a Coulomb friction
element placed in parallel with a viscous damper. In this model,
for nonzero piston velocities $\dot{x}$, the force generated by
the device is given by the fluid as follows:
\begin{eqnarray}\label{Rheo2}
\displaystyle   F(t) = f_c \cdot sgn(\dot{x}) + c_0 \dot{x} + f_0
\end{eqnarray}
where $f_c$ is the coefficient of the frictional force, which is
related to the fluid yield stress, $c_0$ is the damping
coefficient, $f_0$ denoting an offset in the force is included to
account for the nonzero mean observed in the measured force due to
the presence of the accumulator. We remark that if at any point
the velocity of the piston is zero, the force generated in the
frictional element is equal to the applied force. The Bingham
model accounts for electro- and magneto- rheological fluid
behavior beyond the yield point, i.e. for fully developed fluid
flow or sufficiently higher shear rates. However, it assumes that
the fluid remains rigid in the pre-yield region \cite{4}. This the
Bingham model does not describe the fluid's elastic properties at
small deformations and low shear rates which is necessary for
dynamic applications. Lee and Wereley \cite{10} and Wang and
Gordaninejad \cite{11} employed the Herschel-Bulkley model to
accomodate fluid post-yield shear thinning and shear thickening.
In this model, the constant post-yield plastic viscosity in the
Bingham model is replaced with a power law model dependent on the
shear strain rate. However due to its simplicity, the Bingham
model is still very effective, especially in the damper design
phase \cite{9}. The Herschel-Bulkley model can be expressed by
\begin{eqnarray}\label{Rheo3}
\displaystyle   \tau = \left[ \tau_y(H) \cdot sgn(\dot{\gamma}) +
k \left( \dot{\gamma} \right)^{\frac{1}{m}} \right]
sgn(\dot{\gamma})
\end{eqnarray}
where $k$ is the consistency parameter and $m$ is fluid behavior
index of the magnetorheological fluid. For $m > 1 $, Eq.
(\ref{Rheo3}) represents a shear thinning fluid while shear
thickening fluids are described by $m < 1$. For $m=1$, the
Herschel-Bulkley model reduces to the Bingham model \cite{12}.
Zubieta et al. \cite{13} have proposed field-dependent plastic
models, for magnetorheological fluid based on the original Bingham
plastic and Herschel-Bulkley plastic models. In the field
dependent Bingham and Herschel-Bulkley model, the rheological
properties of magnetorheological fluid depend applied magnetic
field and can be estimated by the following equation \cite{12}:
\begin{eqnarray}\label{Rheo4}
\displaystyle   Y = Y_{\infty} + (Y_0 - Y_{\infty}) \left( 2 e^{-
B \alpha_{_{YS}}} - e^{- 2 B \alpha_{_{YS}}} \right)
\end{eqnarray}
where $Y$ stands for a rheological parameters of
magnetorheological fluid such as yield stress, post-yield
viscosity, fluid consistency and flow behavior index. The value of
parameter $Y$ tends from the zero applied field value $Y_0$ to the
saturation value $Y_{\infty}$ and $\alpha_{_{YS}}$ is the
saturation moment index of the $Y$ parameter, $B$ is the applied
magnetic density. The value of $Y_0$, $Y_{\infty}$ and
$\alpha_{_{YS}}$ are
determined from experimental results using curve fitting method.\\
\smallskip
\begin{figure}[!h]
\centering  \includegraphics[width=3.4in]{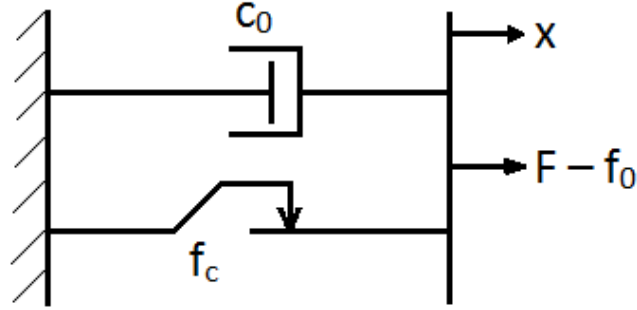}\\
  \caption{Bingham model} \label{fig:1}
\end{figure}

The Bingham body model, presented in Fig. 2 differs from the
Bingham model (Fig. 1) by the introducing of a spring $k$. The
Bingham body model contains in parallel three elements that is
connecting the elements of St. Venant, Newton and the element of
Hooke. To a certain value of the applied force $f_c$ - static
friction force of the St. Venant element, only the spring will
deform, similarly to the elastic Hooke body.
\smallskip
\begin{figure}[!h]
\centering    \includegraphics[width=3.4in]{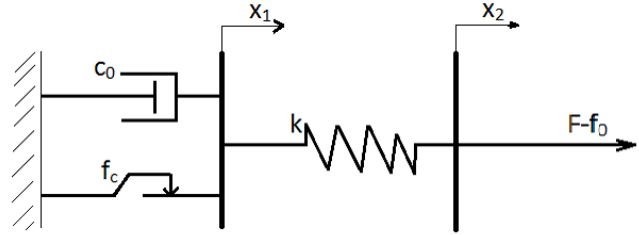}\\
  \caption{Bingham body model} \label{fig:2}
\end{figure}
If this force is greater than $f_c$ the Bingham body will
elongate. The rate of the deformation will be proportional to the
difference of the applied force and the friction force of the St.
Venant element \cite{14}. In theis case the damping force $F$ can
be expressed as
\begin{eqnarray}\label{Rheo5}
\displaystyle   F(t) = \left \{ 
\begin{array}{lll}
f_c sgn (\dot{x}_1) + c_0 \dot{x}_1 + f_0 & \textrm{for} & |F| > f_c \\
k(x_2-x_1)+f_0 & \textrm{for} & |F| \leq f_c\\
\end{array}
\right.
\end{eqnarray}
where $k$ is the stiffness of the elastic body (Hooke model) and
the other parameters have the same meaning as in Eq.
(\ref{Rheo2}).

An extension of the Bingham model is formulated by Gamota and
Filisko \cite{15}. This extension describes the electrorheological
fluid behavior in the pre-yield and post-yield region as well as
the yield point. This viscoelastic-plastic model depends on
connection of the Bingham, Kelvin-Voight body and Hooke body
models (Fig. 3).\\
\smallskip
\begin{figure}[!h]
\centering    \includegraphics[width=3.4in]{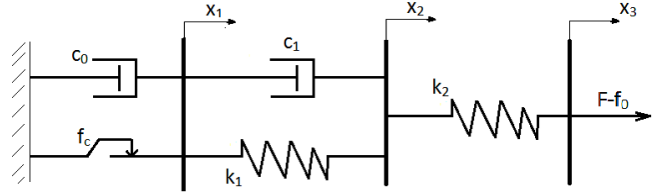}\\
  \caption{Extended Bingham model for the Gamota-Filisko model}
  \label{fig:3}
\end{figure}\\

The damping force in the Gamota-Filisko model is given by
\begin{eqnarray}\label{Rheo6}
\displaystyle   F(t) = \left \{ 
\begin{array}{llll}
k_1(x_2-x_1)+c_1(\dot{x}_2-\dot{x}_1) +f_0  & =  c_0 \dot{x}_1 +  f_c sgn (\dot{x}_1) + f_0 =   \\
                                            & =   k_2(x_3-x_2) + f_0  \ \  \textrm{for} \ |F| > f_c \\
k(x_2-x_1)+ c_1 \dot{x}_2 + f_0 & = k_2(x_3-x_2) + f_0 \ \ \textrm{for} \ |F| \leq f_c \\
\end{array}
\right.
\end{eqnarray}
where $c_0$, $f_0$ and $f_c$ are known from Bingham model
(\ref{Rheo2}) and the parameters $c_1$, $k_1$ and $k_2$ are
associated with the fluid's elastic properties in the pre-yield
region. We remark that if $|F| \leq f_c$ then $\dot{x}_1=0$, which
means that when the friction force $f_c$ related with the new
stress in the fluid is greater than the damping force $F$, the
piston remain motionless. Another view of visco-elastic-plastic
properties of MR damper behavior is proposed by Li et al.
\cite{16}. In essence the damping force is equal to the
visco-plastic force, to which besides the friction force connected
with the fluid shear stress $f_c$, the viscotic force and inertial
force contribute, which can be written in the form \cite{14}
\begin{eqnarray}\label{Rheo7}
\displaystyle   F(t) = f_c \cdot sgn(\dot{x}) + c_0 \dot{x} + m
\ddot{x}
\end{eqnarray}
where $c_0$ is a co-factor of viscotic friction and $m$ is the
mass of replaced MR fluid dependent on the amplitude and frequency
of a kinematic excitation applied to the piston.

A discrete element model with similar components, referred to as
the BingMax model, is reviewed by Makris et al. \cite{17}. It
consists of a Maxwell element in parallel with a Coulomb friction
element as depicted in Fig 4.\\
\smallskip
\begin{figure}[!h]
\centering    \includegraphics[width=3.4in]{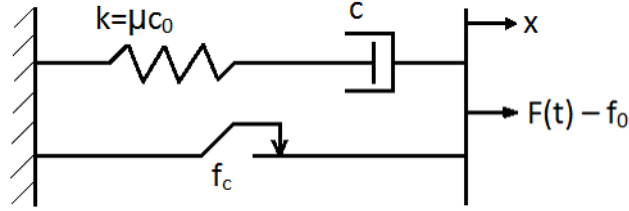}\\
  \caption{BingMax model} \label{fig:4}
\end{figure}

The force $F(t)$ is given by
\begin{eqnarray}\label{Rheo8}
\displaystyle   F(t) = c \cdot \int_0^t \mu e^{-\mu(t-\tau)}
\dot{x}(\tau)\ d\tau + f_c sgn (\dot{x}(t))
\end{eqnarray}
where $\mu$ is a parameter meaning the non-viscous damping effect
(or the relaxation parameter \cite{18}). The model for the damping
force expressed as
\begin{eqnarray}\label{Rheo9}
\displaystyle   f_d(t) = c \cdot \int_0^t \mu e^{-\mu(t-\tau)}
\dot{x}(\tau)\ d\tau
\end{eqnarray}
was originally proposed by Biot \cite{19} and later used by
several authors in the context of dynamics of viscoelastic
systems. Eq. (\ref{Rheo9}) physically implies that the previous
time histories of the velocity $\dot{x}$ contribute to the current
damping force and the most recent instances of velocity have the
highest influence. In the limiting case when $\mu \rightarrow
\infty$, the exponential kernel function approaches the Dirac
delta function $\delta(t)$. For this special case damping force
given by (\ref{Rheo9}) reduces to the case of viscous damping. For
viscoelastic systems, an equation similar to Eq. (\ref{Rheo9}) is
often associated with the stiffness parameter. Cremer and Heckl
\cite{20} concluded that: "of the many after-effect functions,
that are possible in principle, one - the so-called relaxation
function - is physically meaningful".

Based on the above considerations, in what follows we consider the
behavior of a magnetorheological damper after Bingham models which
in addition contains a nonlinear term.

The motion equation is established in the form:
\begin{eqnarray}\label{Rheo10}
\displaystyle   m \frac{d^2 x(\bar{t})}{d\bar{t}^2} + \bar{c}
\cdot \int_0^{\bar{t}} \bar{\mu}
e^{-\bar{\mu}(\bar{t}-\bar{\tau})} \dot{x}(\bar{\tau})\
d\bar{\tau} + k x(\bar{t}) + \alpha k x^3(\bar{t}) + \nonumber \\
+ k \beta sgn(\dot{x}(\bar{t})) + k f_0 =0
\end{eqnarray}
where $m$, $\bar{c}$, $\bar{\mu}$, $k$, $\alpha k$, $\beta k$ and
$k f_0$ are mass, damping relaxation, linear stiffness, nonlinear
stiffness, coefficient of frictional force and offset of the
force, respectively; $x$, $\dot{x}$ and $\ddot{x}$ are the dynamic
responses of the structure (displacement, velocity and
acceleration).

The initial conditions are:
\begin{eqnarray}\label{Rheo11}
\displaystyle   x(0)=A, \ \ \ \dot{x}(0) = v_0 \cdot
\sqrt{\frac{k}{m}}
\end{eqnarray}

The objective of the present paper is to propose an accurate
procedure to nonlinear differential equation of the nonlinear
Bingham model given by Eq. (\ref{Rheo10}), using OHAM. A version
of the OHAM is applied in this study to derive highly accurate
analytical expressions of the solutions using only one iteration
and a small number of steps. Our procedure is independent of the
presence of any small or large parameters, contradistinguishing
from other known methods in literature. The main advantage of this
approach is the control of the convergence of approximate
solutions in a very rigorous way. A very good agreement was found
between our approximate solutions and numerical results, which
proves that our method is very efficient and accurate.

\section{Basic ideas of the optimal homotopy asymptotic method}
\label{sec:2}

\noindent \par Eq. (\ref{Rheo10}) with initial conditions
(\ref{Rheo11}) can be written in a more general form
\begin{eqnarray}\label{Rheo12}
\displaystyle    N\left[x(t)\right]=0
\end{eqnarray}
where $N$ is a given nonlinear differential operator depending on
the unknown function $x(t)$, subject to the initial conditions
\begin{eqnarray}\label{Rheo13}
\displaystyle    B\left(x(t), \ \frac{d x(t)}{d t}\right) = 0.
\end{eqnarray}
Let $x_0(t)$ be an initial approximation of $x(t)$ and $L$ an
arbitrary linear operator such as
\begin{eqnarray}\label{Rheo14}
\displaystyle    L\left[x_0(t)\right] =  0, \ \ B\left(x_0(t), \
\frac{d x_0(t)}{dt}\right) =  0.
\end{eqnarray}
It should be emphasize that this linear operator $L$ is not
unique.

If $p \in [0, \ 1]$ denotes an embedding parameter and $X$ is an
analytic function, then we propose to construct a homotopy
\cite{21} - \cite{26}:
\begin{eqnarray}\label{Rheo15}
\displaystyle   \mathcal{H}\left[L\left(X(t, p)\right), \ H(t,
C_i), \ N\left(X(t,p)\right)\right], \ \ i=1,2,...,s
\end{eqnarray}
with properties
\begin{eqnarray}\label{Rheo16}
\displaystyle   \mathcal{H}\left[L\left(X(t,0)\right), \ H(t,
C_i), \ N\left(X(t,0)\right)\right]= L\left(X(t,0)\right) =
L\left(x_0(t)\right)=0
\end{eqnarray}
\begin{eqnarray}\label{Rheo17}
\displaystyle    \mathcal{H}\left[L\left(X(t,1)\right), \ H(t,
C_i), \ N\left(X(t,1)\right)\right]= H(t, C_i)
N\left(x(t)\right) = 0, \\
i=1,2,...,s \nonumber
\end{eqnarray}
where $H(t, C_i) \neq 0$, is an arbitrary auxiliary
convergence-control function depending on variable $t$ and on $s$
arbitrary parameters $C_1$, $C_2$, ..., $C_s$ unknown now and will
be determined later.

Let us consider the function $X$ in the form
\begin{eqnarray}\label{Rheo18}
\displaystyle    X(t,p) = x_0(t) + p x_1(t, C_i).
\end{eqnarray}

By substituting Eq. (\ref{Rheo18}) into equation obtained by means
of homotopy (\ref{Rheo15})
\begin{eqnarray}\label{Rheo19}
\displaystyle   \mathcal{H}\left[L\left(X(t,p)\right), \ H(t,
C_i), \ N\left(X(t,p)\right)\right]= 0, \ \ i=1, ..., s
\end{eqnarray}
and then equating the coefficients of $p^0$ and $p^1$, we obtain
\begin{eqnarray}\label{Rheo20}
\displaystyle   \mathcal{H}\left[L\left(X(t,p)\right), \ H(t,
C_i),
\ N\left(X(t,p)\right)\right]= \nonumber \\
\displaystyle  = L\left(x_0(t)\right) + p
\left[L\left(x_1(t,C_i)\right) - L\left(x_0(t)\right) -
H(t, C_i)  N\left(x_0(t)\right) \right]=0, \\
i=1,2, ..., s. \nonumber
\end{eqnarray}

From the last equation, we obtain the governing equation of
$x_0(t)$ given by Eq. (\ref{Rheo14}) and the governing equation of
$x_1(t, C_i)$, i.e.:
\begin{eqnarray}\label{Rheo21}
\displaystyle   L\left(x_1(t,C_i)\right) = H(t, C_i)
N\left(x_0(t)\right), \ \ \ B\left(x_1(t,C_i), \ \frac{d
x_1(t,C_i)}{dt}\right) =  0, \\
i=1, \ ..., \ s \nonumber
\end{eqnarray}
where we find the following expression for the nonlinear operator:
\begin{eqnarray}\label{Rheo22}
\displaystyle   N\left(x_0(t)\right) = \sum_{i=1}^{m} h_i(t)g_i(t)
\end{eqnarray}

In the Eq. (\ref{Rheo22}) the functions $h_i(t)$ and $g_i(t)$,
$i=1, \ ..., \ m$ are known and depend on the initial
approximation $x_0(t)$ and also on the nonlinear operator, $m$
being a known integer number.

In this way, taking into account Eq. (\ref{Rheo17}), from Eq.
(\ref{Rheo18}) for $p=1$, we obtain the first-order approximate
solution which becomes\\
\begin{eqnarray}\label{Rheo23}
\displaystyle   \overline{x}(t, C_i) = x_0(t) +  x_1(t, C_i), \ \
\  i=1, \ ..., \ s
\end{eqnarray}

It should be emphasized that $x_0(t)$ and $x_1(t, C_i)$ are
governed by the linear Eqs. (\ref{Rheo14}) and (\ref{Rheo21})
respectively with initial / boundary conditions that come from the
original problem. It is known that the general solution of
nonhomogeneous linear Eq. (\ref{Rheo21}) is equal to the sum of
general solution of the corresponding homogeneous equation and of
some particular solutions of the nonhomogeneous equation. However,
the particular solutions are readily selected only in the
exceptional cases.

In what follows we do not solve Eq. (\ref{Rheo21}), but from the
theory of differential equations, taking into considerations the
method of variation of parameters, Cauchy method, method of
influence function, the operator method and so on, is more
convenient to consider the unknown function $x_1(t, C_i)$, in
the form\\
\begin{eqnarray}\label{Rheo24}
\displaystyle    x_1(t, C_j) = \sum_{i=1}^{n} H_i(t, h_j(t),
C_j)g_i(t),  \ \ \ j=1, \ ..., \ s\\
\displaystyle  B\left(x_1(t,C_i), \ \frac{d x_1(t,C_i)}{dt}\right)
=  0 \nonumber
\end{eqnarray}
where within expression of $H_i(t, h_j(t), C_j)$ appear
combinations of some functions $h_j$, the some terms which are
given by the corresponding homogeneous equation and the unknown
parameters $C_j$, $j=1, \ ..., \ s$. In the sum given by Eq.
(\ref{Rheo24}) appear an arbitrary number of $n$ of the such
terms.

We have a large freedom to choose the value of $n$. We cannot
demand that $x_1(t, C_i)$ to be solutions of Eq. (\ref{Rheo21})
but $\overline{x}(t, C_i)$ given by Eq. (\ref{Rheo23}) with
$x_1(t, C_i)$ given by Eq. (\ref{Rheo24}), are the solutions of
the Eq. (\ref{Rheo12}). This is underlying idea of our method. The
convergence of the approximate solution $\overline{x}(t, C_i)$
given by Eq. (\ref{Rheo23}) depends upon the auxiliary functions
$H_i(t,h_i,C_j)$, $j=1, \ ..., \ s$. There are many possibilities
to choose these functions $H_i$. We try to choose the auxiliary
functions $H_i$ so that within Eq. (\ref{Rheo24}) the term
$\sum_{i=1}^{n} H_i(t,h_j(t),C_j) g_i(t)$ be of the same shape
with the term $\sum_{i=1}^{m} h_i(t) g_i(t)$ given by Eq.
(\ref{Rheo22}). The first-order approximate solution
$\overline{x}(t, C_i)$ also depend on the parameters $C_j$, $j=1,
\ ..., \ s$. The values of these parameters can be optimally
identified via various methods, such as: the least-square method,
the Galerkin method the collocation method, the Ritz method or
minimizing the square residual error:\\
\begin{eqnarray}\label{Rheo25}
\displaystyle    J(C_1, C_2, ..., C_s) = \int_{D}R^2(t, C_1, C_2,
..., C_s) \ dt
\end{eqnarray}
where the residual $R$ is given by\\
\begin{eqnarray}\label{Rheo26}
\displaystyle    R(t, C_1, C_2, ..., C_s) = N\left(\overline{x}(t,
C_i)\right).
\end{eqnarray}

The unknown parameters $C_1, \ C_2, \ ..., \ C_s$ can be
identified from the conditions:\\
\begin{eqnarray}\label{Rheo27}
\displaystyle   \frac{\partial J}{\partial C_1} \ = \
\frac{\partial J}{\partial C_2} \ = \ ... \ = \ \frac{\partial
J}{\partial C_s} = \ 0.
\end{eqnarray}

With these parameters known (called optimal convergence-control
parameters), the first-order approximate solution given by Eq.
(\ref{Rheo23}) is well-known.

It should be emphasized that our procedure contains the auxiliary
functions $H_i(t,f_i,C_j)$, $i=1, \ ..., \ m$, $j=1, \ ..., \ s$
which provides us with a simple way to adjust and control the
convergence of the approximate solutions. It is very important to
properly choose these functions $H_i(t,f_i,C_j)$ which appear in
the construction in the first-order approximation.

\section{Application of OHAM to nonlinear Bingham fluid dampers}
\label{sec:3}

\noindent \par In what follows, we apply our procedure to obtain
approximate solution of Eqs. (\ref{Rheo10}) and (\ref{Rheo11}).
For this purpose, we introduce the dimensionless variables
\begin{eqnarray}\label{Rheo28}
\displaystyle  t= \bar{t} \sqrt{\frac{m}{k}}, \ \ \ \tau=
\bar{\tau} \sqrt{\frac{m}{k}}
\end{eqnarray}
such that Eq. (\ref{Rheo10}) can be expressed as
\begin{eqnarray}\label{Rheo29}
\displaystyle   \ddot{x}(t) + 2 c \cdot \int_0^t \mu
e^{-\mu(t-\tau)} \dot{x}(\tau) \ d\tau + x(t) + \alpha x^3(t) +
\beta sgn(\dot{x}(t)) + f_0 =0
\end{eqnarray}
where $c= \frac{\bar{c}}{2 \sqrt{km}}$ $\mu = \sqrt{\frac{m}{k}}
\bar{\mu}$ and the overdot represents differentiation with respect
to dimensionless time. The initial conditions (\ref{Rheo11})
become
\begin{eqnarray}\label{Rheo30}
\displaystyle   x(0)=A, \ \ \ \dot{x}(0) = v_0
\end{eqnarray}

Making the transformation
\begin{eqnarray}\label{Rheo31}
\displaystyle   x(t)=A e^{-\lambda t} z(t)
\end{eqnarray}
where $\lambda$ is an unknown parameter at this moment, Eq.
(\ref{Rheo29}) can be written as
\begin{eqnarray}\label{Rheo32}
\displaystyle   \ddot{z} -2 \lambda \dot{z} + (\lambda^2 + 2 \mu c
+1)z - 2 c \mu^2 \cdot \int_0^t e^{(\lambda-\mu)(t-\tau)} z(\tau)
\
d\tau +  \nonumber \\
+ \alpha A^2 e^{-2 \lambda t} z^3 + \frac{\beta}{A} e^{\lambda t}
sgn(\dot{z} - \lambda z) + \frac{f_0 e^{\lambda t}}{A} - 2c \mu
e^{(\lambda - \mu ) t}=0
\end{eqnarray}

The initial conditions (\ref{Rheo30}) become
\begin{eqnarray}\label{Rheo33}
\displaystyle   z(0)=1, \ \ \ \dot{z}(0) = \lambda + \frac{v_0}{A}
\end{eqnarray}

For the nonlinear differential equation (\ref{Rheo32}), we choose
the linear operator of the form:\\
\begin{eqnarray}\label{Rheo34}
\displaystyle   L\left(z(t)\right) = \ddot{z}(t) + \omega^2 z(t)
\end{eqnarray}
where $\omega$ is unknown parameter.

The nonlinear operator corresponding to Eqs. (\ref{Rheo32}) and
(\ref{Rheo34}) is\\
\begin{eqnarray}\label{Rheo35}
\displaystyle   N\left(z(t)\right) = -2 \lambda \dot{z} +
(\lambda^2 + 2 \mu c +1 - \omega^2)z - 2 c \mu^2 \cdot \int_0^t
e^{(\lambda-\mu)(t-\tau)} z(\tau) \
d\tau  + \nonumber \\
+ \alpha A^2 e^{-2 \lambda t} z^3 + \frac{\beta}{A} e^{\lambda t}
sgn(\dot{z} - \lambda z) + \frac{f_0 }{A} e^{\lambda t} - 2c \mu
e^{(\lambda - \mu ) t} \qquad
\end{eqnarray}

The initial approximation $z_0(t)$ can be obtained from Eq.
(\ref{Rheo14}) with initial / boundary conditions:\\
\begin{eqnarray}\label{Rheo36}
\displaystyle   z_0(0)=1, \ \ \ \dot{z}_0(0) = -
\frac{\omega^2}{\lambda}
\end{eqnarray}

Eq. (\ref{Rheo14}) with the linear operator (\ref{Rheo34}) and
with initial / boundary conditions (\ref{Rheo36}) has the solution
\begin{eqnarray}\label{Rheo37}
\displaystyle   z_0(t) = \cos \omega t - \frac{\omega}{\lambda}
\sin \omega t
\end{eqnarray}

The nonlinear operator corresponding to Eq. (\ref{Rheo22}) it
holds that
\begin{eqnarray}\label{Rheo38}
\displaystyle   N\left(z_0(t)\right) =  \frac{\lambda (1 +
\lambda^2 + \omega^2 + 2\mu c) - 2\mu ^2 c (\omega^2 - \lambda^2
+\lambda \mu)}{\lambda (\omega^2+(\mu -\lambda)^2)} \cos \omega t
+ \nonumber \\
+ \frac{\omega (\lambda^2+\omega^2-2\mu c
-1)(\lambda^2+\omega^2+\lambda^2-2\mu \lambda)+ 2 \mu^2 c \omega
(\mu -2 \lambda)}{\lambda
(\omega^2+(\mu -\lambda)^2)} \sin \omega t + \nonumber \\
+ \frac{2 \mu^2 c (\omega^2 - \lambda^2) - 2\mu  c \lambda
(\omega^2 + \lambda^2 -2 \lambda \mu)}{\lambda
(\omega^2+(\mu -\lambda)^2)} e^{(\lambda - \mu)t} + \frac{f_0}{A} e^{\lambda t} + \nonumber \\
+ \frac{\alpha A^2}{4} e^{-2 \lambda t} \left[
\frac{3(\omega^2+\lambda^2)}{\lambda^2} \cos \omega t  + \frac{3
\omega (\omega^2+\lambda^2)}{\lambda^3} \sin \omega t +
\frac{\lambda^2 - 3 \omega^2}{\lambda^2} \cos 3\omega t  + \right. \nonumber \\
\left. + \frac{\omega (\omega^2-3\lambda^2)}{\lambda^3} \sin 3
\omega t \right] - \frac{4 \beta}{\lambda A} e^{\lambda t} \left(
\cos \omega t - \frac{1}{3} \cos 3\omega t + \frac{1}{5} \cos
5\omega t - \right. \nonumber \\
\left. - \frac{1}{7} \cos 7\omega t + \frac{1}{9} \cos 9\omega t -
... \right) \qquad
\end{eqnarray}

From Eqs. (\ref{Rheo38}), (\ref{Rheo33}) and (\ref{Rheo24}), we
choose
\begin{eqnarray}\label{Rheo39}
\displaystyle   z_1(t) = C_1 t^2  \cos \omega t + (C_2
t^2 + C_3 t) \sin \omega t + \nonumber \\
+ \left[ \left( C_4 t^2 + (v_0 +
\frac{\omega^2+\lambda^2}{\lambda}A + 3\lambda C_7 ) t \right)
\cos 3\omega t + \right. \nonumber \\
\left. + (C_5 t^2 +C_6 t) \sin 3\omega t +C_7 \right] e^{-2
\lambda t} - C_7 e^{\lambda t} \cos \omega t
\end{eqnarray}

The first-order approximate solution of Eqs. (\ref{Rheo29}) and
(\ref{Rheo30}) given by Eq. (\ref{Rheo23}) is obtained from Eqs.
(\ref{Rheo37}), (\ref{Rheo39}) and (\ref{Rheo31}):\\
\begin{eqnarray}\label{Rheo40}
\displaystyle   x(t) = \left[ (C_1 t^2 + A) \cos \omega t + (C_2
t^2 + C_3 t - \frac{\omega}{\lambda} A) \sin \omega t  \right] e^{- \lambda t} + \nonumber \\
+ \left\{ \left[ C_4 t^2 + \left( \Big(v_0 +
\frac{\omega^2+\lambda^2}{\lambda} \Big)A + 3\lambda C_7 \right) t
\right] \cos 3\omega t + \right. \nonumber \\
\left. + (C_5 t^2 +C_6 t) \sin 3\omega t +C_7 \right\} e^{-3
\lambda t} - C_7 \cos \omega t
\end{eqnarray}

\section{Numerical results}
\label{sec:4}

\noindent \par We illustrate the accuracy of our procedure for the
following values of the parameters: $c=0.1$, $\mu=20$, $\alpha=1$,
$\beta=0.1$, $f_0=0.1$ $A=5$, $v_0=0.1$. The optimal
convergence-control parameters are determined by means of the
least-square method in the three steps as follows:

For $t \in [0, \ 7/2]$ we obtain\\

$$C_1 = -8.3835528344, \ C_2 = 0.4059363155, \ C_3 =
9.8101433224,$$ $$C_4 = 12.7300955924, \ C_5 = -8.3640772984, \
C_6 = -14.6926248464,$$ $$C_7 = -15.2775231617, \ \lambda =
0.4221369200, \ \omega = 1.1700000000.$$

The first-order approximate solution given by Eq. (\ref{Rheo40})
becomes for this first step:\\
\begin{eqnarray}\label{Rheo41}
\displaystyle   \bar{x}(t) = \left[ (- 8.3835528344 t^2 + 5) \cos
\omega t + ( 0.4059363155 t^2 + 9.8101433224 t - \right. \nonumber \\
\left. - 13.8580629241) \sin \omega t \right] e^{-0.4221369200 t}
+ \left[ (12.7300955924 t^2 - \right. \nonumber \\
 - 1.0230014969 t) \cos 3 \omega t + ( - 8.3640772984 t^2
 - 14.692624846 t ) \sin 3 \omega t - \nonumber \\
\left. -15.2775231617 \right]e^{-1.2664107600 t} + 15.2775231617
\cos \omega t \qquad
\end{eqnarray}

For the second step, when $t \in [7/2, \ 7]$ we obtain\\
\begin{eqnarray}\label{Rheo42}
\displaystyle   \bar{x}(t) = \left[ (- 18.7159298670 t^2 +
1.3931962652) \cos \omega t + ( 2.5809325241 t^2 -  \right. \nonumber \\
\left. - 12.3307330274 t - 7.2827523763) \sin \omega t \right]
e^{-0.3726878974 t}
+ \left[ (16.1687136688 t^2 - \right. \nonumber \\
 - 76.6284415838 t) \cos 3 \omega t + ( 1.8324463202 t^2
 - 15.4843220143 t ) \sin 3 \omega t - \nonumber \\
\left. - 68.7655829044 \right]e^{- 1.1180636924 t} + 68.7655829044
\cos \omega t \qquad
\end{eqnarray}
where $\omega = 1.9481775386$, $\lambda = 0.3726878974$.
\smallskip
\begin{figure}[!t]
\centering  \includegraphics[width=3.4in]{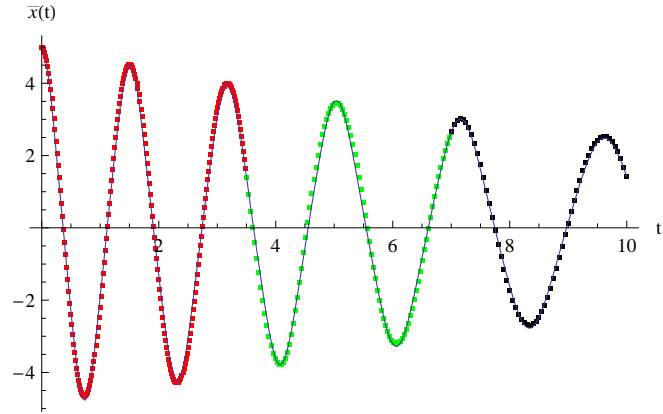}\\
 \caption{Comparison between the approximate solution (\ref{Rheo41}), (\ref{Rheo42}) and (\ref{Rheo43}) and numerical
 solution:  ----- numerical solution;  ......... approximate solution}
 \label{fig:5}
\end{figure}

In the last case, when $t \in [7, \ 10]$, the first-order
approximate solution is\\
\begin{eqnarray}\label{Rheo43}
\displaystyle   \bar{x}(t) = \left[ (0.1521594849 t^2 +
2.6425574994) \cos \omega t + ( - 1.4223903751 t^2 +  \right. \nonumber \\
\left. + 13.8451447750 t - 10.3560770110) \sin \omega t \right]
e^{- 0.1870115667 t}
+ \left[ (- 7.9954565850 t^2 + \right. \nonumber \\
 + 12.3657070213 t) \cos 3 \omega t + ( 1.5073646864 t^2
 - 0.3337501284 t ) \sin 3 \omega t + \nonumber \\
\left. + 1.1037184321 \right]e^{- 0.5610347000 t} - 1.1037184321
\cos \omega t \qquad
\end{eqnarray}
where $\omega = 0.7328908403$, $\lambda = 0.1870115667$.\\

In Fig. 5 is plotted a comparison between the first-order
approximate solution and numerical results.

It can be seen that the solution obtained by the proposed
procedure are nearly identical with the numerical solution
obtained using a fourth-order Runge-Kutta method.

\section{Conclusions}
\label{sec:5}

\noindent \par The Optimal Homotopy Asymptotic Method is employed
to propose an analytic approximate solutions for the nonlinear
Bingham model. Our procedure is valid even if the nonlinear
differential equation does not contain any small or large
parameters. In construction of the homotopy appear some
distinctive concepts as: the linear operator, the nonlinear
operator, the auxiliary functions $H_i(t, C_i)$ and several
optimal convergence-control parameters $\lambda$, $\omega$, $C_1$,
$C_2$, ... which ensure a fast convergence of the solutions. The
example presented in this work, leads to the conclusion that the
obtained results are of very accurate using only one iteration and
three steps. The OHAM provides us with a simple and rigorous way
to control and adjust the convergence of the solution through the
auxiliary functions $H_i(t, C_i)$ involving several parameters
$\lambda$, $\omega$, $C_1$, $C_2$, ... which are optimally
determined. Actually, the capital strength of OHAM is its fast
convergence, which proves that our procedure is very efficient in
practice.

\textbf{Conflict of Interests} \\
The authors declare that there is
no conflict of interests regarding the publication of this paper.

\end{document}